\documentclass[12pt,twoside,doublespace]{article}
\usepackage{latexsym,amssymb}

\def\smskip{\par\vskip 5 pt}
\def\QED{\hfill $\Box$\smskip}

\newtheorem{theorem}{Theorem}[section]
\newtheorem{lemma}{Lemma}[section]

\newtheorem{proposition}{Proposition}[section]

\begin{document}

\begin{center}

\vspace{15pt}

{\Large \bf The Method of Pairwise Variations}

\vspace{10pt}

{\Large \bf with Tolerances for Linearly Constrained}

\vspace{10pt}

{\Large \bf Optimization Problems}

\vspace{35pt}

{\sc I.V.~Konnov\footnote{\normalsize Department of System Analysis
and Information Technologies, Kazan Federal University, ul.
Kremlevskaya, 18, Kazan 420008, Russia. }}

\end{center}

\begin{abstract}
We consider a method of pairwise variations for smooth optimization
problems, which involve polyhedral constraints.
It consists in making steps with respect to the difference of two selected
extreme points of the feasible set together with special threshold control and tolerances
whose values reduce sequentially. The method is simpler and more flexible
than the well-known conditional gradient method, but keeps
its useful sparsity properties and is very suitable for large dimensional
optimization problems. We establish its convergence
under rather mild assumptions. Efficiency of the method is confirmed by
its convergence rates and results of computational  experiments.

{\bf Key words:} Optimization problems; polyhedral feasible set;
pairwise variations; conditional gradient method; threshold control.

\end{abstract}

%11111111111111111111111111111111111111111111111111111111111111111111111111

\section{Introduction} \label{s1}

The usual optimization problem consists in finding the minimal
value of some goal function $f : \mathbb{R}^{m} \to \mathbb{R}$ on
a feasible set $D$ such that $D \subseteq  \mathbb{R}^{m}$. For brevity, we write
this problem as
\begin{equation} \label{eq:1.1}
 \min \limits _{x \in D} \to f(x),
\end{equation}
its solution set is denoted by $D^{*}$ and the optimal value of the
function by $f^{*}$, i.e.
$$
f^{*} = \inf \limits _{x \in D} f(x).
$$
We shall consider a special class of optimization problems, where the set $D$ is a nonempty polyhedron
 and the function $f$ is supposed to be smooth on $D$,
i.e., it is bounded and defined by affine constraints, e.g.
$$
D=\left\{ x\in \mathbb{R}^{m} \ \vrule \ \langle q^{i},x \rangle \leq \beta_{i}, \ i=1, \ldots, l \right\},
$$
where $\langle q,x\rangle$ denotes the usual scalar product of $q$ and $x$.
Then problem (\ref{eq:1.1}) has a solution.

The conditional gradient method is one of the oldest methods, which can be applied for
the above problem. It was first suggested in
\cite{FW56} for the case when the goal function is quadratic and further was developed by many authors;
see e.g. \cite{LP66,DR68,PD78,Dun80}. We recall that the main idea of this method consists in
linearization of the goal function. That is, given the current iterate $x^{k}\in D$, one finds some
solution $y^{k}$ of the problem
\begin{equation} \label{eq:1.2}
\min_{y \in D} \to \langle f'(x^{k}),y \rangle
\end{equation}
and defines $p^{k}=y^{k}-x^{k}$ as a descent direction at $x^{k}$.
Taking a suitable stepsize $\lambda_{k} \in (0,1]$, one sets $x^{k+1}=x^{k}+\lambda_{k}p^{k}$
and so on.

During rather long time, this method was not considered as very efficient because of its
relatively slow convergence in comparison with Newton and projection type methods.
However, it became very popular recently due to several
features significant for many applications, where huge dimensionality and inexact data
create certain drawbacks for more rapid methods.
In particular, its auxiliary linearized problems of form (\ref{eq:1.2})
appear simpler essentially than
the quadratic ones of the most other methods. Next, it
usually yields so-called sparse approximations of a solution with few non-zero
components; see e.g. \cite{Cla10,Jag13}.
 Many efforts were directed to enhance the convergence properties of
the conditional gradient method; see e.g. \cite{GM86,BT04,Jag13,FG16} and the references therein.
In particular, inserting the so-called away steps enabled one to attain the linear rate
of convergence for some classes of optimization problems significant for applications;
see e.g. \cite{GM86,NFS14,LJ15,BS16}.

In this paper, we intend to present some other modification of the conditional gradient method,
which seems more flexible and reduces the total computational expenses.
The main idea follows the bi-coordinate descent method with special threshold control
and tolerances for optimization problems with simplex constraints that
was proposed in \cite{Kon16b}. Unlike the previous methods, its direction choice
requirements are relaxed essentially, which admits different implementation versions.
Its more detailed comparison with the other methods is given in Section \ref{s5}.

In the next section, we give several basic properties of problem (\ref{eq:1.1}),
which will be used for the substantiation of the method.
In Section \ref{s3}, we describe the new method and prove its convergence
in the general case. In Section \ref{s4}, we specialize its convergence properties
for the case where the gradient of the goal function is
Lipschitz continuous, propose some simplifications and
obtain the complexity estimate of the method. In Section \ref{s5},
we discuss its implementation  issues and provide its
comparison with the previously known methods. Section \ref{s6} describes the
results of computational experiments.

%22222222222222222222222222222222222222222222222222222222222222222222222222

\section{Preliminary properties}\label{s2}

We start our consideration from recalling the well known
optimality condition; see e.g. \cite[Theorem 11.1]{Kon07}.

%============================lm:2.1======================================

\begin{lemma} \label{lm:2.1}

(a)  Each solution of problem (\ref{eq:1.1}) is a solution of the
variational inequality (VI for short): Find a point $x^{*} \in D$
such that
\begin{equation} \label{eq:2.1}
\langle f'(x^{*}),x-x^{*} \rangle \geq 0 \quad \forall x \in D.
\end{equation}

(b)  If $f$ is convex, then  each solution of VI
(\ref{eq:2.1}) solves  problem (\ref{eq:1.1}).
\end{lemma}

We denote by $D^{0}$ the solution set of VI
(\ref{eq:2.1}), its elements are called stationary points of problem (\ref{eq:1.1}).

We intend to specialize optimality conditions for problem (\ref{eq:1.1}).
First we note that
\begin{equation} \label{eq:2.1a}
D=\left\{ x\in \mathbb{R}^{m} \ \vrule \ x=\sum \limits_{i\in I} u_{i}z^{i}, \ \sum \limits_{i\in I} u_{i}=1,
\ u_{i}\geq 0, \ i \in I \right\},
\end{equation}
where $z^{i}$ is the $i$-th extreme point (vertex) of the  polyhedron $D$,
$I$ is the set of indices of its extreme points, which is finite, i.e., we can set
$I=\{ 1, \ldots, n \}$.
Given a point $x \in D$, we can hence define the corresponding vector
of weights $u(x)=(u_{1}(x), \dots,u_{n}(x))^{\top}$
of some its associated representation
\begin{equation} \label{eq:2.2}
 x=\sum \limits_{i\in I} u_{i}(x)z^{i}, \ \sum \limits_{i\in I} u_{i}(x)=1,
\ u_{i}(x)\geq 0, \ i \in I.
\end{equation}
Clearly, $u(x)$ is not defined uniquely in general.
Now we give the useful property of solutions of linear programming (LP for short) problems;
see \cite[Section 3.3]{YG69}.

%============================lm:2.2======================================

\begin{lemma} \label{lm:2.2} Let $c$ be a fixed vector in $\mathbb{R}^{m}$.

(i) If a point $x^{*}$ is a solution of
the LP problem
\begin{equation} \label{eq:2.3}
 \min \limits _{x \in D}  \to  \langle c,x \rangle ,
\end{equation}
and
\begin{equation} \label{eq:2.4}
 x^{*}=\sum \limits_{i\in I} u^{*}_{i}z^{i}, \ \sum \limits_{i\in I} u^{*}_{i}=1,
\ u^{*}_{i}\geq 0, \ i \in I;
\end{equation}
then
\begin{equation} \label{eq:2.5}
 \langle c,z^{i} \rangle \left\{ {
\begin{array}{ll}
\displaystyle
\geq \langle c, x^{*} \rangle \quad & \mbox{if} \ u^{*}_{i}=0, \\
=\langle c, x^{*} \rangle  \quad & \mbox{if} \ u^{*}_{i}>0,
\end{array}
} \right. \quad \mbox{for} \ i \in I.
\end{equation}

 (ii) If a point $x^{*} \in D$ satisfies conditions (\ref{eq:2.5})
 for some representation (\ref{eq:2.4}), then it solves problem (\ref{eq:2.3}).
\end{lemma}
{\bf Proof.} Let a point $x^{*}$ be a solution of problem (\ref{eq:2.3})
and (\ref{eq:2.4}) holds.  By definition,
\begin{equation} \label{eq:2.6}
 \langle c, x^{*} \rangle \leq \langle c, z^{i} \rangle, \ \forall  i \in I.
\end{equation}
Define the index sets
 $I_{+}= \{ i \in I \ | \ u^{*}_{i}>0 \}$ and $I_{0}= \{ i \in I \ | \ u^{*}_{i}=0 \}$ and choose
$s \in I_{+}$. Then $u^{*}_{s}>0$ and
\begin{eqnarray*}
\langle c, x^{*} \rangle &=&\sum \limits_{i\in I_{+}} u^{*}_{i}\langle c, z^{i} \rangle
   = u^{*}_{s}\langle c, z^{s} \rangle + \sum \limits_{i\in I_{+}, i \neq s} u_{i}\langle c, z^{i} \rangle \\
 &\geq & u^{*}_{s}\langle c, z^{s} \rangle +(1-u^{*}_{s})\langle c, x^{*} \rangle.
\end{eqnarray*}
It follows that
$\langle c, x^{*} \rangle \geq \langle c, z^{s} \rangle$,
hence
$\langle c, x^{*} \rangle = \langle c, z^{s} \rangle$
in view of (\ref{eq:2.6}). Assertion (i) is true.

Conversely, let a point $x^{*} \in D$ satisfy conditions (\ref{eq:2.5})
for some representation (\ref{eq:2.4}). Take an arbitrary point $x \in D$ and
some associated weight vector $v=u(x)$, then
$$
x=\sum \limits_{i\in I} v_{i}z^{i}, \ \sum \limits_{i\in I} v_{i}=1,
\ v_{i}\geq 0, \ i \in I.
$$
It follows from (\ref{eq:2.5}) that
$$
\langle c, x \rangle=\sum \limits_{i\in I} v_{i}\langle c, z^{i} \rangle \geq \langle c, x^{*} \rangle \sum \limits_{i\in I} v_{i}
=\langle c, x^{*} \rangle,
$$
and assertion (ii) holds true.
\QED

Now we are ready to give optimality conditions for VI (\ref{eq:2.1}), hence for
problem (\ref{eq:1.1}).

%============================pro:2.1======================================

\begin{proposition} \label{pro:2.1}

A point $x^{*}$ with representation (\ref{eq:2.4})
is a solution of VI (\ref{eq:2.1}) if and only if it satisfies each
of the following equivalent conditions:
\begin{eqnarray}
&& x^{*} \in D,  \ \langle f'(x^{*}), z^{i} \rangle \left\{ {
\begin{array}{ll}
\displaystyle
\geq \langle f'(x^{*}), x^{*} \rangle \quad & \mbox{if} \ u^{*}_{i}=0, \\
=\langle f'(x^{*}), x^{*} \rangle \quad & \mbox{if} \ u^{*}_{i}>0,
\end{array}
} \right. \quad \mbox{for} \ i \in I; \label{eq:2.7}\\
&& x^{*} \in D, \ \forall i,j \in I, \ \langle f'(x^{*}), z^{i} \rangle > \langle f'(x^{*}), z^{j} \rangle  \ \Longrightarrow
\ u^{*}_{i}=0; \label{eq:2.8}\\
&& x^{*} \in D, \ \forall i,j \in I, \ u^{*}_{i}>0 \ \Longrightarrow
\ \langle f'(x^{*}), z^{i} \rangle \leq \langle f'(x^{*}), z^{j} \rangle. \label{eq:2.9}
\end{eqnarray}
\end{proposition}
{\bf Proof.} From Lemma \ref{lm:2.2} we clearly have that
a point $x^{*}$ with some representation (\ref{eq:2.4})
is a solution of VI (\ref{eq:2.1}) if and only if it satisfies (\ref{eq:2.7}). Clearly,
(\ref{eq:2.7}) implies (\ref{eq:2.8}) and (\ref{eq:2.8}) implies (\ref{eq:2.9}).
 Let now a  point $x^{*} \in D$ with $u^{*}=u(x^{*})$ satisfy
(\ref{eq:2.9}). Then there exists an index $k$ such that
$u^{*}_{k}>0$. Set
$$
 \alpha=\min_{i \in I} \langle f'(x^{*}), z^{i} \rangle.
$$
Then (\ref{eq:2.9}) implies $\langle f'(x^{*}), z^{i} \rangle  = \alpha$ if
$u^{*}_{i}>0$ and $\langle f'(x^{*}), z^{i} \rangle  \geq \alpha$ if $u^{*}_{i}=0$,
hence (\ref{eq:2.7}) holds. \QED

Given a number $\varepsilon>0$ and a point $x \in D$ with some associated weight
vector $u(x)$ from (\ref{eq:2.2}), let
$$
  I_{\varepsilon}(x)=\{ i \in I \ | \
     u_{i}(x) \geq \varepsilon \}.
$$
The number of vertices may be too large, however, we can evaluate the weights
implicitly from feasible step-sizes.

%============================pro:2.2======================================

\begin{proposition} \label{pro:2.2}  Given $x \in D$ and $i \in I$, let
\begin{equation} \label{eq:2.10}
 x+\alpha(z^{j}-z^{i}) \notin D  \  \mbox{for some} \ j \in I, \ \mbox{if} \ \alpha \geq \varepsilon >0.
\end{equation}
Then $u_{i}(x) < \varepsilon$ for any weight
vector $u(x)$.
\end{proposition}
{\bf Proof.} On the contrary, suppose that (\ref{eq:2.10}) holds, but there exists
a weight vector $u=u(x)$ with $u_{i} \geq \varepsilon$. Take $\alpha=u_{i}$
and arbitrary $j \in I$. Then can define the point $y= x+\alpha(z^{j}-z^{i})$ such that
$$
y=\sum \limits_{s\in I} u_{s}z^{s}+\alpha(z^{j}-z^{i})= \sum \limits_{s\in I} v_{s}z^{s},
$$
where
$$
v_{s}= \left\{ {
\begin{array}{ll}
\displaystyle
0, \quad & \mbox{if} \ s=i, \\
u_{j}+u_{i}, \quad & \mbox{if} \ s=j, \\
u_{s}, \quad & \mbox{otherwise};
\end{array}
} \right.
$$
besides,
$$
\sum \limits_{s\in I} v_{s}=1, \ v_{s}\geq 0, \ s \in I.
$$
It follows that $y \in D$, a contradiction. \QED

%33333333333333333333333333333333333333333333333333333333333333333333333333

\section{Method and its convergence}\label{s3}

The method of pairwise variations with tolerances (PVM for short)
for VI (\ref{eq:2.1}) is described as follows.
Let $\mathbb{Z}_{+}$ denote the set of non-negative integers.

\medskip \noindent
%=================Method of pairwise variations==================================
 {\bf Method (PVM).} \\ {\em Initialization:} Choose a point $w^{0} \in D$, numbers $\beta \in (0,1)$,
$\theta  \in (0,1)$, and
sequences $\{\delta _{l}\} \searrow 0$, $\{\varepsilon _{l}\}
\searrow 0$ with $\varepsilon _{0} \in (0,1)$. Set $l:=1$.\\
{\em Step 0:} Set  $k:=0$, $x^{0}:=w^{l-1}$.\\
{\em Step 1:}  Choose an index $i \in I_{\varepsilon _{l}}(x^{k})$
for some associated weight vector $u^{k}=u(x^{k})$ and an index $ j \in
I$ such that
\begin{equation} \label{eq:3.2}
\langle f'(x^{k}), z^{i}-z^{j} \rangle  \geq \delta _{l},
\end{equation}
choose $\gamma_{k}\in [\varepsilon _{l},u^{k}_{i}]$, set $i_{k}:=i$, $j_{k}:=j$ and go to Step 2. Otherwise
(i.e. if (\ref{eq:3.2}) does not hold for
all $i \in I_{\varepsilon _{l}}(x^{k})$ associated to some weight vector $u(x^{k})$  and $ j \in
I$) set $w^{l}:=x^{k}$, $l:=l+1$ and go to Step 0. {\em (Restart)} \\
{\em Step 2:}  Set $d^{k}:= z^{j_{k}}-z^{i_{k}}$,
determine $m_{k}$ as the smallest number in $\mathbb{Z}_{+}$ such that
\begin{equation} \label{eq:3.3}
 f (x^{k}+\theta ^{m_{k}}\gamma_{k} d^{k})
 \leq f (x^{k})+\beta \theta ^{m_{k}}\gamma_{k}
  \langle f'(x^{k}),d^{k} \rangle,
\end{equation}
set $\lambda_{k}:=\theta ^{m_{k}}\gamma_{k}$, $x^{k+1}:=x^{k}+\lambda_{k}d^{k}$,
$k:=k+1$ and go to Step 1.\\
\medskip

Thus, the method has a two-level structure where each outer iteration (stage) $l$
contains some number of inner iterations in $k$
with the fixed tolerances $\delta _{l}$ and  $\varepsilon _{l}$. Completing each stage,
that is marked as restart, leads to decrease of their values.

Note that $i_{k}\neq j_{k} $ due to (\ref{eq:3.2}), besides, $\gamma_{k} \geq \varepsilon_{l}$
and the point $x^{k}+\gamma_{k} d^{k}$ is always feasible.  Moreover, by definition,
\begin{equation} \label{eq:3.5}
\mu _{k} =\langle f'(x^{k}),d^{k}
\rangle=\langle f'(x^{k}), z^{j_{k}}-z^{i_{k}} \rangle \leq
 -\delta_{l}<0,
\end{equation}
in (\ref{eq:3.3}). It follows that
 \begin{equation} \label{eq:3.6}
 f (x^{k+1}) \leq f (x^{k})+\beta \lambda_{k} \mu _{k} \leq f (x^{k})-\beta \lambda_{k}\delta_{l},
\end{equation}

We now justify the linesearch.

%============================lm 3.1=====================================================

\begin{lemma} \label{lm:3.1} The linesearch procedure in Step 2 is always finite.
\end{lemma}
{\bf Proof.}
If we suppose that the linesearch procedure is infinite, then (\ref{eq:3.3}) does not hold and
$$
(\theta ^{m_{k}}\gamma_{k})^{-1 }(f (x^{k}+\theta ^{m_{k}}\gamma_{k}
d^{k}) - f (x^{k}))>\beta \mu _{k},
$$
for $m_{k} \to \infty$. Hence, by taking the limit we have $    \mu _{k}
\geq \beta \mu _{k}$, hence $\mu _{k} \geq 0$, a
contradiction with $\mu _{k}  \leq -\delta_{l}<0$ in (\ref{eq:3.5}). \QED

We show that each stage is well defined.

%=========================pro:3.1==================================================

\begin{proposition} \label{pro:3.1}
The number of iterations at each stage $l$ is finite.
\end{proposition}
{\bf Proof.}
Fix any $l$. Since the sequence $\{x^{k}\}$ is bounded, it has limit points.
Besides, by (\ref{eq:3.6}), we have $f^{*}\leq f(x^{k})$ and $f(x^{k+1})\leq f(x^{k})-\beta \delta_{l}
\lambda_{k}$, hence
$$
\lim \limits_{k\rightarrow \infty }\lambda_{k}=0.
$$
Suppose that the sequence $\{x^{k}\}$ is infinite. Since the set $I$
is finite, there is a pair of indices $(i_{k},j_{k})=(i,j)$,
 which is repeated infinitely. Take the corresponding
subsequence $\{k_{s}\}$, then $d^{k_{s}}=\bar d=z^{j}-z^{i}$.
Without loss of generality, we can
suppose that the subsequence $\{x^{k_{s}}\}$ converges to a point
$\bar x$ and due to (\ref{eq:3.5}) we have
$$
\langle f'(\bar x),\bar d \rangle= \lim \limits_{s\rightarrow \infty
} \langle f'(x^{k_{s}}),\bar d \rangle \leq -\delta_{l}.
$$
However,  (\ref{eq:3.3}) does not hold for the step-size $\lambda_{k}/\theta$.
Setting $k=k_{s}$ gives
$$
(\lambda_{k_{s}}/\theta)^{-1 }(f
(x^{k_{s}}+(\lambda_{k_{s}}/\theta) \bar d) - f
(x^{k_{s}}))>\beta \langle f'(x^{k_{s}}),\bar d \rangle,
$$
hence, by taking the limit $s\rightarrow \infty$ we obtain
$$
\langle f'(\bar x),\bar d \rangle= \lim \limits_{s\rightarrow \infty
} \left\{(\lambda_{k_{s}}/\theta)^{-1 }(f
(x^{k_{s}}+(\lambda_{k_{s}}/\theta) \bar d) - f
(x^{k_{s}}))\right\} \geq  \beta \langle f'(\bar x),\bar d
\rangle,
$$
i.e.,  $   (1-\beta ) \langle f'(\bar x),\bar d \rangle \geq 0$,
which is a contradiction. \QED

We are ready to prove convergence of the whole method.

%3===================thm 3.1=========================================

\begin{theorem} \label{thm:3.1}
Under the assumptions made it holds that:

(i)  the number of changes of index $k$ at each stage $l$  is finite;

(ii) the sequence $\{w^{l}\}$ generated by method (PVM) has limit points, all
these limit points are solutions of VI (\ref{eq:2.1});

(iii) if $f$ is convex, then
\begin{equation} \label{eq:3.7}
 \lim \limits_{l\rightarrow \infty} f(w^{l})=f^{*};
\end{equation}
and all the limit points of $\{w^{l}\}$ belong to $D^{*}$.
\end{theorem}
{\bf Proof.} Assertion (i) has been obtained in Proposition \ref{pro:3.1}.
By construction, the sequence $\{w^{l}\}$ is bounded, hence it has limit points.
Moreover, $f(w^{l+1})\leq f(w^{l})$, hence
\begin{equation} \label{eq:3.8}
\lim \limits_{l\rightarrow \infty }f(w^{l})=\mu.
\end{equation}
Take an arbitrary limit point $\bar w$ of $\{w^{l}\}$, then $\bar w\in D$,
$$
\lim \limits_{t\rightarrow \infty }w^{l_{t}}=\bar w.
$$
By definition (\ref{eq:2.1a}), each point $w^{l}$ is associated with some weight
vector $v^{l}=u(w^{l})$  such that
$$
w^{l} =\sum \limits_{s\in I} v^{l}_{s}z^{s}, \ \sum \limits_{i\in I} v^{l}_{s}=1, \  v^{l}_{s} \geq 0, \ s \in I.
$$
Clearly, the sequence $\{v^{l}\}$ is bounded and must have limit points.
Without loss of generality we can suppose that
$$
\bar v =\lim \limits_{t\rightarrow \infty }v^{l_{t}},
$$
then
$$
\bar w =\sum \limits_{s\in I} \bar v_{s}z^{s}, \ \sum \limits_{i\in I} \bar v_{s}=1, \  \bar v_{s} \geq 0, \ s \in I.
$$
For $l>0$ we must have
$$
  \langle f'(w^{l}), z^{i}-z^{j} \rangle \leq \delta_{l}
  \ \mbox{for all} \ i,j \in I \  \mbox{with} \ v^{l}_{i} \geq \varepsilon_{l}.
$$

Let $p$ be an arbitrary index such that $\bar v_{p} >0$. Then
$v^{l_{t}}_{p} \geq \varepsilon_{l_{t}}$ for $t$ large enough, hence
$$
 \langle f'(w^{l_{t}}), z^{p}-z^{q} \rangle \leq \delta_{l_{t}} \ \mbox{for all}
  \ q \in I.
$$
Taking the limit $t\rightarrow \infty$, we obtain
$$
\langle f'(\bar w),  z^{p}-z^{q} \rangle \leq 0 \ \mbox{for all}
  \ q \in I.
$$
This means that the point $\bar w$ satisfies the optimality
conditions (\ref{eq:2.9}). Due to Proposition \ref{pro:2.1}, $\bar
w$ solves VI (\ref{eq:2.1}) and assertion (ii) holds. Next, if $f$ is
convex, then by Lemma \ref{lm:2.1} each limit point
of $\{w^{l}\}$ belongs to $D^{0}$, hence $\mu=f^{*}$ in (\ref{eq:3.8}).
This gives (\ref{eq:3.7}) and assertion (iii).
\QED

%444444444444444444444444444444444444444444444444444444444444444444444444444

\section{Convergence in the Lipschitz gradient case}\label{s4}

The above descent method is very flexible and admits various
modifications and extensions. In particular, we can take the
exact one-dimensional minimization rule instead of the current
Armijo rule in (\ref{eq:3.3}). The convergence then  can be obtained
along the same lines; see e.g. \cite[Section 6.1]{Kon13}.

If the gradient of the function $f$ is Lipschitz continuous on $D$ with
some constant $L >0$, i.e., $\| f'(y)-f'(x)\| \leq L \|y-x \|$
for any vectors $x$ and $y$,  we can take the useful property of
such functions
$$
f (y) \leq f(x)+\langle
f'(x),y-x \rangle +0.5L\|y-x \|^{2};
$$
see \cite[Chapter III, Lemma 1.2]{DR68}.
This gives us an explicit lower bound for the step-size. In fact, at Step 2 we have
$$
 f (x^{k}+\lambda d^{k}) - f (x^{k}) \leq \lambda [\langle f'(x^{k}),d^{k} \rangle
    +0.5L \lambda\|d^{k}\|^{2} ]  \leq \beta \lambda \langle f'(x^{k}),d^{k} \rangle,
$$
if $ \lambda \leq -(1-\beta)\langle f'(x^{k}),d^{k} \rangle /(L\|d^{k}\|^{2})$.
However, $\langle f'(x^{k}),d^{k} \rangle \leq -\delta_{l}$ at stage $l$,
besides, $ \|d^{k}\| \leq B= {\rm Diam} D < \infty$.
If we take $ \lambda_{k} =\lambda \delta_{l}$ with $ \lambda \in (0, \bar \lambda]$ and
$$
 \bar \lambda =\min\{(1-\beta)/(LB^{2}), \varepsilon_{l}\},
$$
 then
\begin{equation} \label{eq:4.1}
 f (x^{k}+\lambda_{k} d^{k})
 \leq f (x^{k})+\beta \lambda_{k} \langle f'(x^{k}),d^{k} \rangle,
\end{equation}
as desired. In such a way we can drop the line-search procedure in
Step 2.  Obviously, the assertions of
Proposition \ref{pro:3.1} and Theorem
\ref{thm:3.1} remain true for this version.
This version reduces the computational expenses essentially but require the evaluation of the
Lipschitz constants. We can use several approaches to avoid this drawback.

 Firstly, we can apply the step-size rule $ \lambda_{k} =\varepsilon_{l} \delta_{l}$ at
 stage $l$ without any line-search. Then $\varepsilon_{l} \leq \bar \lambda$ for $l$ large enough and
 the convergence can be proved as in the previous case since the values of the function $f$
 are bounded from above on the compact set $D$. Besides, after the finite number of stages we
 will have the basic inequality $f(w^{l+1})\leq f(w^{l})$, which implies (\ref{eq:3.8}).

 Secondly, we can apply the divergent step-size rule
\begin{equation} \label{eq:4.1a}
 \sum \limits_{k=0}^{\infty }\lambda_{k}=\infty, \
 \sum \limits_{k=0}^{\infty }\lambda^{2}_{k}<\infty, \ \lambda_{k} \in (0, \varepsilon_{l}], \ k=1,2, \ldots,
\end{equation}
at stage $l$. For instance, we can set $\lambda_{k} = \varepsilon_{l}/(k+1)$. Then again
$\lambda_{k} \leq \bar \lambda \delta_{l}$ for $k$ large enough.
Then the assertion of Proposition \ref{pro:3.1}  remains true. In fact, if we suppose
that the sequence $\{x^{k}\}$ is infinite, (\ref{eq:3.6}) gives
$f(x^{s})\leq f(x^{s-1})-\beta \delta_{l}
\lambda_{s}$, hence
$$
f^{*}\leq f(x^{k}) \leq f(x^{0})-\beta \delta_{l} \sum \limits_{s=0}^{k }\lambda_{s},
$$
which is a contradiction. Then assertion (ii) of Theorem
\ref{thm:3.1} can be proved as above. Assertion (iii) follows from (ii)
and the continuity of $f$. Therefore, rule (\ref{eq:4.1a}) also provides convergence.

Of course, there is no necessity now to evaluate the Lipschitz constant and diameter of $D$.

Due to Lemma \ref{lm:2.1}, the value
$$
\Delta (x)=\max _{y\in D}\langle f'(x),x-y \rangle
$$
gives a gap function for VI (\ref{eq:2.1}). We intend to
obtain an error bound for VI (\ref{eq:2.1}) at $w^{l}$.
Since $D$ is compact, we can define
$$
 \sigma=\max _{i \in I}\max _{x\in D} \langle f'(x),z^{i} \rangle.
$$

%============================pro:4.1======================================

\begin{proposition} \label{pro:4.1} For each stage $l$, we have
\begin{equation} \label{eq:4.2}
 \Delta (w^{l}) \leq  \delta_{l} +2 n \varepsilon_{l} \sigma.
\end{equation}
\end{proposition}
{\bf Proof.} By definition,
$$
\min _{y\in D}\langle f'(w^{l}),y \rangle =\langle f'(x),z^{t} \rangle
$$
for some $t \in I$. We recall that $u(w^{l})=v^{l}$ and $I_{\varepsilon_{l}}(w^{l})=\{ i \in I \ | \
     v^{l}_{i} \geq \varepsilon_{l} \}$. It follows that
\begin{eqnarray*}
\Delta (w^{l}) &=& \sum_{s \in I}v^{l}_{s}\langle f'(w^{l}),z^{s} \rangle -\langle f'(x),z^{t} \rangle
      = \sum_{s \in I} v^{l}_{s}\langle f'(w^{l}),z^{s} - z^{t} \rangle\\
    &=& \sum_{s \in I_{\varepsilon_{l}}(w^{l})}  v^{l}_{s}\langle f'(w^{l}),z^{s} - z^{t} \rangle
         +\sum_{s \notin I_{\varepsilon_{l}}(w^{l})}  v^{l}_{s}\langle f'(w^{l}),z^{s} - z^{t} \rangle\\
       &\leq &  \delta_{l} \sum_{s \in I_{\varepsilon_{l}}(w^{l})}  v^{l}_{s}
          +\varepsilon_{l} \sum_{s \notin I_{\varepsilon_{l}}(w^{l})} \langle f'(w^{l}),z^{s} - z^{t} \rangle \\
       &\leq &   \delta_{l} +2 n \varepsilon_{l} \sigma .
\end{eqnarray*}
Therefore, estimate (\ref{eq:4.2}) holds true. \QED

As the method has a two-level structure with each stage containing a
finite number of inner iterations, it is more suitable
to derive its complexity estimate, which gives the total amount of
work of the method. We now suppose that the  function $f$ is convex
and its gradient is Lipschitz continuous with constant $L$. For simplicity, we take
the above version with  the fixed stepsize $ \lambda^{k} =\bar \lambda \delta_{l}$.

We take the value $\Phi(x)=f (x)-f^{*}$ as an
accuracy measure for our method. More precisely, given a starting
point $ z^{0}$ and a number $\alpha > 0$, we define
the complexity of the method, denoted by $N(\alpha )$,  as the
total number of inner iterations at $l(\alpha )$ stages such that
$l(\alpha )$ is the maximal number $l$ with
$\Phi(z^{l}) \geq \alpha$, hence,
\begin{equation} \label{eq:4.3}
 N (\alpha ) \leq  \sum ^{l(\alpha  )} _{l=1} N_{(l)},
\end{equation}
where $N_{(l)}$ denotes the total number of iterations at  stage
$l$. We proceed to estimate the right-hand side of (\ref{eq:4.3}).
To change the parameters, we apply the rule
\begin{equation} \label{eq:4.4}
 \delta _{l} = \varepsilon _{l} = \nu ^{l}\delta_{0},  l=0,1,\ldots;
  \quad \nu \in (0,1), \delta_{0}>0.
\end{equation}
By (\ref{eq:4.1}), we have
$$
 f (x^{k+1}) \leq f (x^{k})-\beta \bar \lambda \delta^{2} _{l},
$$
hence
\begin{equation} \label{eq:4.5}
 N_{(l)} \leq \Phi(z^{l-1})/(\beta\bar
\lambda \delta^{2} _{l}).
\end{equation}
Under the above assumptions from Proposition \ref{pro:4.1}  we obtain
$$
 \Phi(z^{l}) = f(z^{l}) - f^{*} \leq \Delta (z^{l}) \leq \delta _{l} + 2n \sigma\varepsilon _{l}
             = \delta_{0}C_{1}\nu ^{l},
$$
where $C_{1}=1+ 2n\sigma $. It follows that
$$
\nu^{-l(\alpha  )} \leq \delta_{0} C_{1}/\alpha .
$$
Besides, using (\ref{eq:4.5}) now gives
$$
N_{(l)} \leq C_{1} \delta_{0} \nu ^{l-1} /(\beta\bar \lambda \nu ^{2l} \delta
_{0}^{2}) = C_{1}  L B^{2} /(\beta(1-\beta) \nu ^{l+1} \delta
_{0})=C_{2}\nu^{-l-1},
$$
where $C_{2}=C_{1} L B^{2}/(\beta(1-\beta)\delta_{0})$.

Combining both the inequalities in (\ref{eq:4.3}), we obtain
\begin{eqnarray*}
N (\alpha ) && \leq C_{2}\nu^{-1}\sum ^{l(\alpha  )}
_{l=1}
\nu^{-l}  \leq C_{2} (\nu^{-l(\alpha  )}-1)/(1-\nu)  \\
&& \leq C_{2} (C_{1}/\alpha-1)/(1-\nu).
\end{eqnarray*}
We have established the complexity estimate.

%===========================thm 4.1=====================================

\begin{theorem} \label{thm:4.1} Let the function $f : X \to \mathbb{R}$ be convex
and its gradient be Lipschitz continuous with constant $L$.
 Let a sequence $\{w^{l}\}$ be generated by (PVM) with
  the stepsize rule $ \lambda_{k} =\bar \lambda \delta _{l}$ at stage $l$. 
  If the parameters satisfy conditions (\ref{eq:4.4}),
the method has the complexity estimate
$$
N (\alpha ) \leq C_{2} (C_{1}/\alpha-1)/(1-\nu),
$$
where $C_{1}=1+ 2n\sigma$ and $C_{2}=C_{1} L B^{2}/(\beta(1-\beta)\delta_{0})$.
\end{theorem}

We see that the above estimate corresponds to those of the usual
conditional gradient methods, which solves the
linearized problem (\ref{eq:1.2}) at each iteration; see
\cite{LP66,Dun80}.

%55555555555555555555555555555555555555555555555555555555555555555555555555555

\section{Implementation issues}\label{s5}

In this section, we discuss some questions of implementation  of (PVM)
for different kinds of feasible sets and provide its comparison with the previously
known methods. In fact, implementation  of (PVM) requires some associated
weight vector $u^{k}=u(x^{k})$ for each iteration point $x^{k}$.
This vector is used for finding a suitable index $i \in I_{\varepsilon _{l}}(x^{k})$.
We again note that it suffices to have an arbitrary weight vector of $x^{k}$.
The first way is to choose such a vector at the starting point $w^{0}$
and change it sequentially in conformity with the iteration process.
For the sake of clarity, we give its full description now.

\medskip \noindent
%=================Method of pairwise variations with weight changes==================================
 {\bf (PVM) with explicit weight changes.} \\ {\em Initialization:} Choose a point $w^{0} \in D$ with some associated
weight vector $v^{0}=u(w^{0})$, numbers $\beta \in (0,1)$,
$\theta  \in (0,1)$, and
sequences $\{\delta _{l}\} \searrow 0$, $\{\varepsilon _{l}\}
\searrow 0$ with $\varepsilon _{0} \in (0,1)$. Set $l:=1$.\\
{\em Step 0:} Set  $k:=0$, $x^{0}:=w^{l-1}$, $u(x^{0}):=u^{0}:=v^{l-1}$.\\
{\em Step 1:}  Choose a pair of indices $i \in I_{\varepsilon _{l}}(x^{k})$ and $ j \in
I$ such that
\begin{equation} \label{eq:5.1}
\langle f'(x^{k}), z^{i}-z^{j} \rangle  \geq \delta _{l},
\end{equation}
set $\gamma_{k}:=u^{k}_{i}$, $i_{k}:=i$, $j_{k}:=j$ and go to Step 2. Otherwise
(i.e. if (\ref{eq:5.1}) does not hold for
all $i \in I_{\varepsilon _{l}}(x^{k})$ and $ j \in
I$) set $w^{l}:=x^{k}$, $u(w^{l}):=v^{l}:=u^{k}$, $l:=l+1$ and go to Step 0. {\em (Restart)} \\
{\em Step 2:}  Set $d^{k}:= z^{j_{k}}-z^{i_{k}}$,
determine $m_{k}$ as the smallest number in $\mathbb{Z}_{+}$ such that
$$
 f (x^{k}+\theta ^{m_{k}}\gamma_{k} d^{k})
 \leq f (x^{k})+\beta \theta ^{m_{k}}\gamma_{k}
  \langle f'(x^{k}),d^{k} \rangle,
$$
set $\lambda_{k}:=\theta ^{m_{k}}\gamma_{k}$, $x^{k+1}:=x^{k}+\lambda_{k}d^{k}$,
\begin{equation} \label{eq:5.2}
u_{s}(x^{k+1}):=u^{k+1}_{s}:= \left\{ {
\begin{array}{ll}
\displaystyle
u^{k}_{s}-\lambda_{k} \quad & \mbox{if} \ s=i_{k}, \\
u^{k}_{s}+\lambda_{k} \quad & \mbox{if} \ s=j_{k}, \\
u^{k}_{s} \quad & \mbox{otherwise};
\end{array}
} \right.
\end{equation}
$k:=k+1$ and go to Step 1.\\
\medskip

Observe that
$$
  I_{\varepsilon _{l}}(x^{k})=\{ i \in I \ | \
     u^{k}_{i} \geq \varepsilon _{l} \}
$$
and that we can simply set $\gamma_{k}:=u^{k}_{i}$ in Step 1. Clearly, formula
(\ref{eq:5.2}) gives the weight vector $u(x^{k+1})$ associated to $x^{k+1}$ without
solution of any system of equations. In fact,
\begin{eqnarray*}
x^{k+1} &=& x^{k}+\lambda_{k}d^{k}=\sum \limits_{i\in I} u^{k}_{i}z^{i}+\lambda_{k}(z^{j_{k}}-z^{i_{k}}) \\
 &=& \sum \limits_{i\in I, i \neq i_{k},j_{k}} u^{k}_{i}z^{i}
   +(u^{k}_{i_{k}}-\lambda_{k})z^{i_{k}}+(u^{k}_{j_{k}}+\lambda_{k})z^{j_{k}}=\sum \limits_{i\in I} u^{k+1}_{i}z^{i};
\end{eqnarray*}
in addition, we have
$$
\sum \limits_{i\in I} u^{k+1}_{i}=1 \ \mbox{and} \ u^{k+1}_{i} \geq 0, \ i \in I.
$$
Therefore, each iterate changes only two components of the current weight vector.
Clearly, it suffices to keep only positive components of this vector.
Set
$$
  I_{+}(x^{k})=\{ i \in I \ | \
     u^{k}_{i} >0 \},
$$
then the number of indices in $I_{+}(x^{k})$ is much more smaller than that in $I$.
For instance,  any segment $[a,b]$ in $\mathbb{R}^{m}$ has $2^{m}$ vertices, whereas
any point $x \in [a,b]$ can be represented by $m+1$ vertices, i.e., for this weight vector, set
$I_{+}(x)$ contains $m+1$ items.

The other way to implementation consists in calculation the necessary
weights from the iteration point $x^{k}$. Both the approaches coincide if each
point $x \in D$ has the unique weight vector $u=u(x)$. This is the case for the simplices.
In fact, take
$$
D=\left\{ x\in \mathbb{R}^{m}_{+} \ \vrule \ \langle a,x\rangle=\tau \right\},
$$
$\tau$ is a fixed positive number, $a$ is a fixed vector with positive coordinates,
$\mathbb{R}^{m}_{+}$ denotes the non-negative orthant in $\mathbb{R}^{m}$.
Given $x \in D$, set
$$
\sigma(x)=\sum \limits_{i=1}^{m} x_{i},
$$
then
$$
z^{i}_{s}= \left\{ {
\begin{array}{ll}
\displaystyle
\tau/a_{s} \quad & \mbox{if} \ s=i, \\
0 \quad & \mbox{otherwise};
\end{array}
} \right.
$$
for $i=1, \ldots, n$ and
$$
u_{s}= \left\{ {
\begin{array}{ll}
\displaystyle
x_{s}/\sigma(x) \quad & \mbox{if} \ s=i, \\
0 \quad & \mbox{otherwise}.
\end{array}
} \right.
$$
However, the second way may be useful if the weight vector $u(x)$ is not defined uniquely.
Moreover, we can evaluate the weight implicitly by using Proposition \ref{pro:2.2}.

Next, the current condition (\ref{eq:3.2}) (or (\ref{eq:5.2})) for
selection of the pair of indices $i_{k}$ and $j_{k}$ can
be implemented within various rules. It seems suitable to find
$z^{i_{k}}$ as an approximate solution of problem (\ref{eq:1.2}) and
$j_{k}$ as an approximate solution of the problem
\begin{equation} \label{eq:5.3}
\max_{s \in I_{\varepsilon _{l}}(x^{k})} \to \langle f'(x^{k}),z^{s}\rangle.
\end{equation}
That is, we can make several steps of any algorithm toward the solutions of
(\ref{eq:1.2}) and (\ref{eq:5.3}) for satisfying (\ref{eq:3.2}).
We can even solve (\ref{eq:1.2}) exactly, and then check the indices from
$I_{+}(x^{k})$ sequentially. This procedure does not seem too difficult since
the number of indices in $I_{+}(x^{k})$ is much more smaller than that in $I$.

We should observe that all these implementations of (PVM) are closely related with the so-called
\lq\lq atomic" or weighting representation (\ref{eq:2.1a}) of the feasible set $D$.
The usual conditional gradient method and its version with away steps
can utilize the standard definition of $D$, whereas their \lq\lq pure" 
weighting versions are also rather popular; see e.g.
\cite{Jag13,LJ15,BS16}. It should be also noticed that all the weighting versions of the methods including (PVM)
can be in principle applied to problem (\ref{eq:1.1}) where the feasible set $D$ is represented as
$$
D=\left\{ x \ \vrule \ x=\sum \limits_{i\in I} u_{i}z^{i}, \ \sum \limits_{i\in I} u_{i}=1,
\ u_{i}\geq 0, \ z^{i} \in H, \ i \in I \right\},
$$
where $H$ is some Hilbert space and the index set $I=\{ 1, \ldots, n \}$ is finite. 
This is treated as linear variable transformation, i.e. $x=Tu$ for 
some linear mapping $T : \mathbb{R}^{n} \to H$, which transforms the 
standard simplex in $\mathbb{R}^{n}$ into $D$. 
In turn, the goal function $f(x)$ is also replaced with the function
$\varphi(u)=f(Tu)$. Since the simplex is convex and compact,
convergence of the method can be proved along the same lines.

It was mentioned in Section \ref{s1} that the main drawback of the
usual conditional gradient method is its rather slow convergence.
Incorporating the away steps, whose calculation requires the solution of the auxiliary problem
\begin{equation} \label{eq:5.4}
\max_{s \in I_{+}(x^{k})} \to \langle f'(x^{k}),z^{s}\rangle,
\end{equation}
enables one to attain the linear rate
of convergence for some classes of optimization problems, but the
computational experiments do not reveal this preference; see e.g.
\cite{GM86,NFS14,BS16}. Instead of these steps we can utilize the so-called
pairwise away or swap directions. Namely, let $z^{j}$ and $z^{i}$ be solutions of
problems (\ref{eq:1.2}) and (\ref{eq:5.4}), respectively. Then, we can take
 $d^{k}= z^{j}-z^{i}$ as the descent direction at the $k$-th iteration; see
\cite{NFS14}. It should be noted that the method based on the same
pairwise directions was first suggested in \cite{DS69} for
network equilibrium problems. In \cite{Kor80}, a similar method was suggested
for general smooth optimization problems with simplex type constraints.
These marginal based index choice methods became very popular after
appearance of their big data applications;
see e.g. \cite{Bec14} for more details and references.
It was also mentioned in Section \ref{s1} that (PVM) can be viewed as an extension of
the bi-coordinate descent method (BCV) with special threshold control proposed in \cite{Kon16b}
for optimization problems with simplex constraints.
That is, (PVM)  can be applied for optimization problems with
arbitrary affine constraints due to the utilization of the weight vectors
which is treated as variable transformation. In comparison with the marginal
 swap direction strategy, (PVM) does not insist on solutions of auxiliary problems
of form (\ref{eq:1.2}) and (\ref{eq:5.4}), which enables us to reduce
the computational expenses significantly. Nevertheless, (PVM) maintains
 the useful sparse iteration point property, as all the mentioned conditional gradient methods.

%66666666666666666666666666666666666666666666666666666666666666666666666666666666666666666666666666

\section{Computational experiments}\label{s6}

In order to check the performance of (PVM) we carried
out computational experiments.  We took also
the usual conditional gradient method (CGM),
the marginal-based swap direction descent method (MDM),
with the same  Armijo linesearch and compared them with (PVM).
They were implemented in Delphi with double precision
arithmetic. The main goal was to compare the numbers of iterations (it) and
calculations of partial derivatives of $f$ (calc) for attaining the same
accuracy $\delta'=0.1$.  We took the following accuracy measure:
$$
\Delta_{k}=\max_{y \in D}\langle f'(x^{k}),x^{k}-y \rangle.
$$
We chose $\beta =\theta =0.5$ for the methods, and the rule
$\delta_{l+1}=\nu \delta_{l}$,
$\varepsilon_{l+1}=\nu\varepsilon_{l}$  with $\nu = 0.5$ for (PVM).

We first took the simplex as the feasible set, i.e.,
\begin{equation} \label{eq:6.1}
D=\left\{ x\in \mathbb{R}^{m}_{+} \ \vrule \ \sum \limits_{i=1}^{m} x_{i}=\tau \right\}.
\end{equation}
We took two starting points, namely,  $x'=(\tau/m)e$
where $e$ denote the vector of units in $\mathbb{R}^{m} $, and
$x''=\tau e^{1}$ where $e^{1}$ denote the first coordinate vector in $\mathbb{R}^{m} $.
Also, we set $\tau=10$.

In the first series, we took the quadratic cost function.
We chose $f(x)=\varphi(x)$ where
\begin{equation} \label{eq:6.2}
\varphi (x)= 0.5 \langle Px,x \rangle -\langle q,x \rangle,
\end{equation}
the elements of the matrix $P$ are defined by
\begin{equation} \label{eq:6.3}
p_{ij}= \left\{ {
\begin{array}{rl}
\displaystyle
\sin (i) \cos (j) \quad & \mbox{if} \ i<j, \\
\sin (j) \cos (i) \quad & \mbox{if} \ i>j, \\
\sum \limits_{i=1}^{m} | p_{ij}| +1 \quad & \mbox{if} \ i=j;
\end{array}
} \right.
\end{equation}
and $q=\mathbf{0}$. The results for the starting points $x'$ and
$x''$ are given in Tables \ref{tbl:1} and \ref{tbl:2}, respectively.
\begin{table}
\caption{Starting point $x'$, quadratic cost function} \label{tbl:1}
\begin{center}
\begin{tabular}{|l|c|c|c|c|}
\hline
                     &  (CGM)  &   (MDM)  & (PVM)    \\
\hline
                     & it / calc  &  it / calc & it / calc \\
                     &           &           &       \\
\hline
    $m=5$            & 202 / 1010  &  11 / 55 & 11 / 53 \\
                     &           &        &       \\
\hline
     $m=10$           & at 500 / 5000            &   34 / 340 & 37 / 279 \\
                      &  $\Delta_{k}=0.25$     &         &   \\
\hline
    $m=20$            & at 500 / 10000            &   49 / 980 & 50 / 703 \\
                       &  $\Delta_{k}=0.11$     &         &   \\
\hline
$m=50$           & at 500 / 25000            &   87 / 4350 & 108 / 3574\\
                    &  $\Delta_{k}=0.39$     &         &   \\
\hline
$m=100$          & at 500 / 50000            &   221 / 22100 & 267 / 17594 \\
                &  $\Delta_{k}=0.62$     &         &   \\
\hline
\end{tabular}
\end{center}
\end{table}
\begin{table}
\caption{Starting point $x''$, quadratic cost function} \label{tbl:2}
\begin{center}
\begin{tabular}{|l|c|c|c|c|}
\hline
                     &  (CGM)  &   (MDM)  & (PVM)    \\
\hline
                     & it / calc  &  it / calc & it / calc \\
                     &           &           &       \\
\hline
    $m=5$            & 47 / 235  &  14 / 70 & 17 / 74 \\
                     &           &        &       \\
\hline
     $m=10$           & 194 / 1940            &   37 / 370 & 42 / 307 \\
                      &              &         &   \\
\hline
    $m=20$            & at 500 / 10000            &   124 / 2480 & 124 / 1668 \\
                       &  $\Delta_{k}=0.44$     &         &   \\
\hline
$m=50$           & at 500 / 25000            &   326 / 16300 & 211 / 7046\\
                    &  $\Delta_{k}=1.22$     &         &   \\
\hline
$m=100$          & at 500 / 50000            &    at 500 / 50000 & 399 / 25213 \\
                &  $\Delta_{k}=2.93$     &       $\Delta_{k}=0.31$   &    \\
\hline
\end{tabular}
\end{center}
\end{table}

In the second series, we took the convex cost function
\begin{equation} \label{eq:6.4}
f(x)=\varphi(x)+1/(\langle c,x\rangle+\mu),
\end{equation}
where the function $\varphi$ was defined as above
in (\ref{eq:6.2})--(\ref{eq:6.3}),
the elements of the vector $c$ are defined by
$$
c_{i}= 2+\sin(i) \ \mbox{ for } \ i=1,\ldots, m,
$$
and $\mu=5$.
The results for the starting points $x'$ and
$x''$ are given in Tables \ref{tbl:3} and \ref{tbl:4}, respectively.
\begin{table}
\caption{Starting point $x'$, convex cost function} \label{tbl:3}
\begin{center}
\begin{tabular}{|l|c|c|c|c|}
\hline
                     &  (CGM)  &   (MDM)  & (PVM)    \\
\hline
                     & it / calc  &  it / calc & it / calc \\
                     &           &           &       \\
\hline
    $m=5$            & 203 / 1015  &  11 / 55 & 11 / 53 \\
                     &           &        &       \\
\hline
     $m=10$           & at 500 / 5000            &   34 / 340 & 38 / 287 \\
                      &  $\Delta_{k}=0.21$     &         &   \\
\hline
    $m=20$            & 491 / 9820            &   53 / 1060 & 46 / 666 \\
                       &                  &         &   \\
\hline
$m=50$           & at 500 / 25000            &   83 / 4150 & 107 / 3427\\
                    &  $\Delta_{k}=0.41$     &         &   \\
\hline
$m=100$          & at 500 / 50000            &   211 / 21100 & 267 / 17012 \\
                &  $\Delta_{k}=0.61$     &         &   \\
\hline
\end{tabular}
\end{center}
\end{table}
\begin{table}
\caption{Starting point $x''$, convex cost function} \label{tbl:4}
\begin{center}
\begin{tabular}{|l|c|c|c|c|}
\hline
                     &  (CGM)  &   (MDM)  & (PVM)    \\
\hline
                     & it / calc  &  it / calc & it / calc \\
                     &           &           &       \\
\hline
    $m=5$            & 44 / 220  &  14 / 70 & 15 / 67 \\
                     &           &        &       \\
\hline
     $m=10$           & 198 / 1980            &   37 / 370 & 43 / 312 \\
                      &              &         &   \\
\hline
    $m=20$            & at 500 / 10000            &   114 / 2280 & 138 / 1839 \\
                       &  $\Delta_{k}=0.45$     &         &   \\
\hline
$m=50$           & at 500 / 25000            &   319 / 15950 & 227 / 7354\\
                    &  $\Delta_{k}=1.24$     &         &   \\
\hline
$m=100$          & at 500 / 50000            &    at 500 / 50000 & 405 / 25758 \\
                &  $\Delta_{k}=2.97$     &       $\Delta_{k}=0.33$   &    \\
\hline
\end{tabular}
\end{center}
\end{table}

Next, we took the more general feasible set instead of (\ref{eq:6.1}):
$$
D=\left\{ x\in \mathbb{R}^{m}_{+} \ \vrule \ \sum \limits_{i=1}^{m} a_{i}x_{i}=\tau \right\}.
$$
the elements of the vector $a$ were defined by
$$
a_{i}= 1.5+\sin(i) \ \mbox{ for } \ i=1,\ldots, m,
$$
and fixed $\tau=10$. We took only the starting point $x''=(\tau/a_{1}) e^{1}$.

In the first series, we took the quadratic cost function from (\ref{eq:6.2})--(\ref{eq:6.3}),
the elements of the vector $q$ were defined by
$$
q_{i}= \sin(i)/i \ \mbox{ for } \ i=1,\ldots, m.
$$
 The results are given in Table \ref{tbl:5}.
\begin{table}
\caption{Quadratic cost function} \label{tbl:5}
\begin{center}
\begin{tabular}{|l|c|c|c|c|}
\hline
                     &  (CGM)  &   (MDM)  & (PVM)    \\
\hline
                     & it / calc  &  it / calc & it / calc \\
                     &           &           &       \\
\hline
    $m=5$            & 20 / 100  &  9 / 45 & 11 / 48 \\
                     &           &        &       \\
\hline
     $m=10$           & 82 / 820            &   29 / 290 & 27 / 210 \\
                      &              &         &   \\
\hline
    $m=20$            & 199 / 3980            &   48 / 960 & 49 / 644 \\
                       &       &         &   \\
\hline
$m=50$           & at 500 / 25000            &   101 / 5050 & 119 / 3630\\
                    &  $\Delta_{k}=0.21$     &         &   \\
\hline
$m=100$          & at 500 / 50000            &    203 / 20300 & 286 / 17080 \\
                &  $\Delta_{k}=0.62$     &          &    \\
\hline
\end{tabular}
\end{center}
\end{table}

In the second series, we took the convex cost function from (\ref{eq:6.4})
where the function $\varphi$ was defined as above.
The results for  are given in Table \ref{tbl:6}.
\begin{table}
\caption{Convex cost function} \label{tbl:6}
\begin{center}
\begin{tabular}{|l|c|c|c|c|}
\hline
                     &  (CGM)  &   (MDM)  & (PVM)    \\
\hline
                     & it / calc  &  it / calc & it / calc \\
                     &           &           &       \\
\hline
    $m=5$            & 20 / 100  &  7 / 35 & 11 / 48 \\
                     &           &        &       \\
\hline
     $m=10$           & 79 / 790            &   27 / 270 & 25 / 189 \\
                      &              &         &   \\
\hline
    $m=20$            & 204 / 4080            &   49 / 980 & 51 / 677 \\
                       &       &         &   \\
\hline
$m=50$           & at 500 / 25000            &   100 / 5000 & 117 / 3618 \\
                    &  $\Delta_{k}=0.19$     &         &   \\
\hline
$m=100$          & at 500 / 50000            &    210 / 21000 & 307 / 18468 \\
                &  $\Delta_{k}=0.64$     &         &       \\
\hline
\end{tabular}
\end{center}
\end{table}

In all the cases, (PVM) showed rather rapid convergence, it outperformed (MDM) in the number of
total calculations if $m \geq 10$, besides, (PVM) and (MDM) appeared better essentially than (CGM).

%777777777777777777777777777777777777777777777777777777777777777777777777777777

\section{Conclusions}

We suggested a new class of descent methods for smooth
optimization problems involving general
affine constraints. The method is based on selective
pairwise variations together with some threshold strategy. 
It keeps the convergence properties of the usual gradient
ones together with reduction of the total computational expenses.
Besides, it is suitable for large scale problems.  The preliminary results of
computational tests show rather rapid and stable convergence of the new method
in comparison with the previous conditional gradient type methods.

%%%%%%%%%%%%%%%%%%%%%%%%%%%%%%%%%%%%%%%%%%%%%%%%%%%%%%%%%%%%%%%%%%%%%%%%%%%%%%%%%%%%%%%%%%%%%

\section*{Acknowledgement}

This work was supported by the RFBR grant, project No. 16-01-00109a
and by grant No. 297689 from Academy of Finland.

%@@@@@@@@@@@@@@  Bibliography  @@@@@@@@@@@@@@@@@@@@@@@@@@@@@@@@@@@@@@@@@@@

\end{document}